\documentclass[12pt]{article} 
\usepackage[utf8]{inputenc} 
\usepackage{geometry}
\usepackage{graphicx}
\usepackage{tikz}
\usetikzlibrary{trees,shapes,arrows}

\usepackage{paralist}
\usepackage{subfig}
\usepackage{sidecap}
\usepackage{array}

\usepackage{amsthm}
\usepackage{amsmath}
\usepackage{amsfonts}
\usepackage{mathtools}

\usepackage{float}
\usepackage{hyperref}
\hypersetup{
	colorlinks=false,
	linkbordercolor=gray,
	citebordercolor=gray,
	urlbordercolor=gray}%
\usepackage{cleveref}
\usepackage[textwidth=1.1in,textsize=footnotesize]{todonotes}

\usepackage{sectsty}
\allsectionsfont{\sffamily\mdseries\scshape\centering\nohang}

\usepackage{abstract}

\usepackage{caption}
\usepackage{sidecap}

\allowdisplaybreaks

\usepackage[switch,modulo]{lineno}
\usepackage{amssymb}
\usepackage[affil-it]{authblk}

\newtheoremstyle{mplain}
{\topsep}   % ABOVESPACE
{\topsep}   % BELOWSPACE
{\itshape}  % BODYFONT
{0pt}       % INDENT (empty value is the same as 0pt)
{\scshape} % HEADFONT
{:}         % HEADPUNCT
{5pt plus 1pt minus 1pt} % HEADSPACE
{}          % CUSTOM-HEAD-SPEC

\renewenvironment{proof}[1][\proofname]{{ \scshape #1. }}{\qed}

\theoremstyle{mplain}

\newtheorem{theorem}{Theorem}[]
\crefname{theorem}{theorem}{theorems} 
\Crefname{theorem}{Theorem}{Theorems}

\newtheorem{lemma}[theorem]{Lemma}
\crefname{lemma}{lemma}{lemmas} 
\Crefname{lemma}{Lemma}{Lemmas}

\newtheorem{proposition}[theorem]{Proposition}
\crefname{proposition}{proposition}{propositions} 
\Crefname{proposition}{Proposition}{Propositions}

\newtheorem{cor}[theorem]{Corollary}
\crefname{cor}{corollary}{corollaries} 
\Crefname{cor}{Corollary}{Corollaries}

\newtheorem{definition}[theorem]{Definition}
\crefname{definition}{definition}{definitions} 
\Crefname{definition}{Definition}{Definitions}

\crefname{example}{example}{examples} 
\Crefname{example}{Example}{Examples}

\newtheorem{remark}[theorem]{Remark}
\crefname{remark}{remark}{remarks} 
\Crefname{remark}{Remark}{Remarks}

\newtheorem{conj}[theorem]{Conjecture}
\crefname{conj}{conjecture}{conjectures} 
\Crefname{conj}{Conjecture}{Conjectures}

\crefname{problem}{problem}{problems} 
\Crefname{problem}{Problem}{Problems}

\crefname{claim}{claim}{claims} 
\Crefname{claim}{Claim}{Claims}

\numberwithin{theorem}{section}

\newcommand{\W}{W}
\newcommand{\E}{\mathbb{E}}

\newcommand{\R}{\mathbb{R}}

\renewcommand{\Pr}{\mathbb{P}}
\newcommand{\Op}{O_\mathbb{P}}
\newcommand{\op}{o_\mathbb{P}}
\newcommand{\Omp}{\Omega_\mathbb{P}}
\newcommand{\omp}{\omega_\mathbb{P}}
\newcommand{\Tp}{\Theta_\mathbb{P}}
\newcommand{\eps}{\varepsilon}

\newcommand{\toprob}{\overset{\mathbb{P}}{\to}}

\newcommand{\Exp}{\mathrm{Exp}}
\newcommand{\e}{\mathrm{e}}

\providecommand{\keywords}[1]{\textbf{\textit{Keywords---}} {\small #1}}

	\title{Minimal H-factors and covers}
	\author{Lorenzo Federico}
	\affil{{\small Department of Political Sciences, LUISS Guido Carli \\ \texttt{lfederico@luiss.it} \\Research supported by the Horizon 2020 Framework Programme through the project Media Futures: Data-driven innovation hub for the media value chain (grant number 951962)}}
	\author{Joel Larsson Danielsson}
	\affil{{ \small Department of Statistics, Lund University School of Economics and Management  \\ \texttt{joel.danielsson@stat.lu.se} \\Research partially supported by the Horizon 2020 Framework Programme through the project Random Graph Geometry and Convergence (grant number 639046)}}

\begin{document}
\maketitle

%-----------------------------------------------------------------------------------------------------------------------------------------------------------------------------------------------------------------------------

\begin{abstract}
Given a fixed small graph $H$ and a larger graph $G$, an $H$-factor is a collection of vertex-disjoint subgraphs $H'\subset G$, each isomorphic to $H$, that cover the vertices of $G$.

If $G$ is the complete graph $K_n$ equipped with independent $U(0,1)$ edge weights, what is the lowest total weight of an $H$-factor? This problem has previously been considered for e.g.\ $H=K_2$.

We show that if $H$ contains a cycle, then the minimum weight is sharply concentrated around some $L_n = \Theta(n^{1-1/d^*})$ (where $d^*$ is the maximum $1$-density of any subgraph of $H$).
Some of our results also hold for $H$-covers, where the copies of $H$ are not required to be vertex-disjoint.
\end{abstract}

\section{Introduction}

\subsection{Threshold and minimum weight problems}
Let $K_n$ denote the complete graph on $n$ vertices, equipped with i.i.d. edge weights $\{ X_e \}_{e\in \mathcal{E}(K_n)}$. We will use the terms `weight' and `cost' interchangeably. For now, let the weight distribution be uniform on $[0,1]$ -- it will turn out that e.g.\ $\Exp(1)$ weights will give the same asymptotic behaviour. For details, see \cref{intro:pseudodim}.
For any family $\mathcal{F}$ of subgraphs of $K_n$, there are two closely related problems.
\begin{description}
    \item[Threshold:] For which $p$ is an $F\in \mathcal{F}$ likely to appear in $G_{n,p}$? That is, if we define the random variable $T:=\min_{F\in \mathcal{F}}\max_{e\in \mathcal{E}(F)} X_e$, what is its distribution? Is it sharply concentrated around its expected value?
    \item[Minimum weight:] The minimal weight of an $F\in \mathcal{F}$ is a random variable $W:=\min_{F\in \mathcal{F}}\sum_{e\in \mathcal{E}(F)} X_e$. What is its distribution? Is it sharply concentrated?
\end{description}

 This pair of problems has been studied for many families $\mathcal{F}$, particularly for families where each $F\in \mathcal{F}$ is spanning -- i.e.\ $\mathcal{V}(F)=\mathcal{V}(K_n)$.
Threshold problems are generally more well-studied than the corresponding minimum-weight problems.
It has been observed that for many natural choices of $\mathcal{F}$, the property of $G_{n,p}$ containing some $F\in \mathcal{F}$ exhibits the sharp threshold phenomenon, i.e.\ $T$ is sharply concentrated around its mean. And for these families, this is often true of the minimum weight $W$ as well.

For instance, if $\mathcal{F}$ is the family of spanning trees, then $T$ is the threshold for connectivity in $G_{n,p}$, and $W$ is the minimal cost of a spanning tree. It's well known that $p={\log n}/{n}$~\cite{threshold:connectivity} is the threshold function for connectivity, and $W\toprob \zeta(3)$~\cite{tree-zeta-3}.
Closely related is the case when $\mathcal{F}$ is the family of perfect matchings. Here the threshold is again $p={\log n}/{n}$~\cite{threshold:matching} (in both cases the minimal obstruction is local and is the existence of an isolated vertex) and $W\toprob\zeta(2)$\cite{matching-zeta-2}.
Similarly for Hamilton cycles, the threshold is  $p={(\log n+\log \log n)}/{n}$~\cite{KOMLOS198355,threshold:hamiltoncycle1,threshold:hamiltoncycle2} and $W\toprob 2.04...$ \cite{tsp-204}.

%Generalizing to hypergraphs, one can study thresholds in the hypergraph generalization of the Erdős–Rényi model -- often called the Linial-Meshulam model -- where each $k$-tuple of vertices occur as an edge independently with probability $p$. 
%Here the threshold for the appearance of a spanning $2$-sphere in a $3$-uniform Linial-Meshulam graph has been determined to be $c/\sqrt{n}$ for an explicit $c$ \cite{2-sphere-threshold}, while in the corresponding minimum weight problem the weight is sharply concentrated around some $a_n=\Theta(\sqrt{n})$.

The goal of this paper is to consider the case when $\mathcal{F}$ is the family of either $H$-factors or $H$-covers. An $H$-factor is a collection of vertex-disjoint subgraphs of $K_n$, each isomorphic to $H$, which collectively cover all $n$ vertices. $H$-covers are defined similarly, but the condition that the subgraphs are vertex-disjoint is dropped. While the threshold version of the $H$-factor problem has received much attention (e.g. \cite{JohKahVu08,rucinski-partialfactor}), the minimum-weight version has (as far as we are aware) not yet been studied.
We prove the following, as well as a similar result for partial factors, and weaker results for covers. These can all be found in \cref{mainthm:cover,mainthm:factor}.
\begin{theorem}
\label{thm:intro-main}
Assume $H$ is a fixed graph with at least one cycle,  $d^*>1$ is its maximum $1$-density as defined in \cref{subsection:density}, and $\Op$ is as defined in \cref{subsection:notation}.

Let the random variable $F_H=F_H(n)$ be the minimum weight of an $H$-factor on $K_n$ (equipped with i.i.d. uniform $[0,1]$ or exponential $\Exp(1)$ edge weights).
Then there exists $M=\Theta(n^{1-1/d^*})$ such that $|F_H-M|=\Op(M^{3/4})$, as $n\to \infty$.
\end{theorem}

Although we work with graphs throughout this paper, in principle our proof method should work for hypergraphs as well, under suitable conditions. However, several theorems we cite have only been proven in the graph setting and would need to be adapted to work for hypergraphs.

\subsection{Proof strategy}
\label{subsection:completingpartialfactors}
Our proof follows a significantly different strategy compared to the study of the minimal perfect matching. The condition that the graph $H$ contains a cycle is equivalent to $d^*>1$.
For such $d^*$, note that the minimum weight of an $H$-factor scales like a positive power of $n$. This scaling enables the following divide-and-conquer approach, which is the main novel contribution of this paper.
It is crucial for two parts our proof to work: the upper bound and sharp concentration of $F_H$. 

A large partial $H$-factor $Q$, covering some $n-k$ vertices, can be completed by adding the lowest weight $H$-factor $Q'$ on the remaining $k$ vertices. Any such $Q$ has weight of order at least $n^{1-1/d^*}$, while $Q'$ has weight of order at most $k^{1-1/d^*}$. So if $k\ll n$, we can complete a large partial factor at a relatively small extra cost. Note that for graphs $H$ with $d^*=1$ (such as $H=K_2$) this does not work, since the minimum weight there instead scales like $F_H=\Theta(1)$. 

However, the $Q$ above might have been picked based on the edge weights (for instance, as the lowest-weight such partial factor) so that the weights of $Q$ and $Q'$ are not independent.
To avoid this dependence, we employ a variant of a trick originally due to Walkup~\cite{walkup1979expected} 
in \cref{section:redgreensplit}: split every edge into a green and red edge, and put independent $\Exp(1-t)$ and $\Exp(t)$ on them respectively, for some small $t>0$. This ensures independence, and the minimum of the two weights on such a pair of edges follows the distribution $\Exp(1)$. We can now find a large partial factor on the green edges (at a slightly inflated cost), and complete the factor using a small number of red edges (at a highly inflated cost).

For the upper bound, we use an upper bound on the cost of a partial factor (due to Rucinski~\cite{rucinski-partialfactor}), and recursively apply the red-green split to find larger partial factors on the remaining vertices. 
To show concentration, we study a dual problem: For some $L=L(n)$, how large is the largest partial $H$-factor with weight at most $L$? We use Talagrand's concentration inequality to show that this size is sharply concentrated around a large value, and then the red-green split trick to complete this large partial factor at small additional cost.

\subsection{Structure of the paper}
We begin with some definitions in \cref{section:defs}. In \cref{section:thms}, we state our main results (\cref{mainthm:cover,mainthm:factor}) and one conjecture, and compare this with previous work.
We then provide proofs in \Cref{section:proofs} under the assumption that the edge weights follow an exponential distribution.
In \cref{section:lowerbound,section:lowerbound-cover} we prove the lower bounds of \cref{mainthm:cover,mainthm:factor} respectively.
\Cref{section:redgreensplit} is devoted to the red-green split trick mentioned in \cref{subsection:completingpartialfactors}.
This trick is then used in \cref{section:upperbound,section:concentration}, where we prove the upper bound and sharp concentration respectively of \cref{mainthm:factor}. In \cref{pseudo-dim} we show that the (asymptotic) distribution of the minimum cost of an $H$-cover or -factor is unchanged if the edge weight distribution is changed from exponential to uniform or some other distribution of pseudo-dimension 1.
Finally, in \cref{section:unconcentrated-cover} we discuss some pathological examples that illustrate why the equivalent of \cref{mainthm:factor} cannot hold for covers.

\section{Definitions and notation}
\label{section:defs}
\subsection{Notation}
\label{subsection:notation}
We will use $\toprob$ to denote convergence in probability, and write $X\overset{d}{=}Y$ if the random variables $X$ and $Y$ follow the same distribution.
We will also use both standard and probabilistic big-$O$ notation. For sequences $X_n,Y_n$ of random variables, the notations $X_n=\Op(Y_n)$ and $Y_n=\Omp(X_n)$ are equivalent, and mean that for any $\eps>0$, there exists a $C=C(\eps)$ such that ${\Pr(|X_n|>C |Y_n|) < \eps}$ for all sufficiently large $n$. Let $X_n = \Tp(Y_n)$ denote that $X_n=\Op(Y_n)$ and $X_n=\Omp(Y_n)$.
Similarly, the notations $X_n=\op(Y_n)$,  $Y_n=\omp(X_n)$ and $X_n \ll Y_n$ are equivalent, and mean that $X_n/Y_n\toprob 0$. When both $X_n$ and $Y_n$ are deterministic, these definitions agree with those for standard big-$O$ notation.

For any graph $G$, we will use $\mathcal{V}(G)$ and $\mathcal{E}(G)$ to refer to its vertex set and edge set respectively, while $v_G:=|\mathcal{V}(G)|$ and $e_G:=|\mathcal{E}(G)|$. Since we will also frequently need to refer to Euler's number $\e\approx 2.78$, we will use a different font to avoid confusion: $\e$ instead of $e$. We will also use $\exp(x)$ for the exponential function, and $\Exp(\lambda)$ for the exponential distribution.

\subsection{Density and balanced graphs}
\label{subsection:density}

For every graph $H$ we define its density as $d_H:={e_H}/{(v_H-1)}$. This quantity is sometimes called the $1$-density (referring to the $-1$ in the denominator), but we will refer to it simply as the density. %We will assume throughout the paper that $H$ is such that $d_H>1$.
We call $H$ \emph{strictly balanced} if $d_G<d_H$ for every subgraph $G \subset H$. Furthermore, let $d^*:=\max\{d_G: G\subseteq H\}$, and let $H^*\subseteq H$ be a subgraph which achieves this maximal density $d^*$.

\subsection{Covers and factors}
\begin{definition}
An $(\alpha,H)$-cover $Q$ is a collection of subgraphs\footnote{Subgraphs, not induced subgraphs.} of the complete graph $K_n$, each of which is isomorphic to $H$, and such that at most $\alpha n$ vertices of $K_n$ are not covered by any copy $H'$ of $H$, i.e.
\[\Big|\!\bigcup_{H'\in Q}\mathcal{V}(H')\Big| \geq (1-\alpha) n.\]
An $(\alpha,H)$-factor is an $(\alpha,H)$-cover such that $\mathcal{V}(H')$ and $\mathcal{V}(H'')$ are disjoint for any two $H',H''\in Q$ with $H'\neq H''$.
For $\alpha=0$ we will refer to $(0,H)$-covers and $(0,H)$-factors simply as $H$-covers and $H$-factors respectively. For $\alpha>0$, we will also refer to $(\alpha,H)$-covers and $(\alpha,H)$-factors as \emph{partial} covers and factors.
\end{definition}
Note that by definition an $H$-factor over $n$ vertices exists if and only if $v_H$ divides $n$, and that for all the valid $H$-factors, $|Q|=n/v_H$.
From now on, we tacitly assume all results about factors to hold only when $v_H$ divides $n$.

\subsection{Edge weight distribution}
\label{intro:pseudodim}
It turns out that the precise distribution of the (positive) edge weights doesn't matter, only its asymptotic behavior near $0$. That is, our results will hold under the following condition: if $F$ is the common cdf of the edge weights, and $F(x)=\lambda x+o(x)$ for some $\lambda>0$ as $x\to 0$. This property is sometimes referred to as $F$ having \emph{pseudo-dimension} 1. For distributions without atoms, this corresponds to a density function tending to $\lambda$ near $0$. Some examples of such distributions are Uniform $U(0,1)$, Exponential, and (for certain values of their parameters) Gamma, Beta and Chi-squared.

We will prove this later in \cref{pseudo-dim}, but for the sake of convenience we will until then assume that the edge weights follow an exponential distribution $\Exp(1)$.

\subsection{Minimum weight covers and factors }

We will also (with minor abuse of notation) let $\mathcal{E}(Q):=\bigcup_{H'\in Q} \mathcal{E}(H')$ denote the (multi-) set of edges that occur in some copy of $H$. For factors this is a set, while for covers this is a multiset where the multiplicity of an edge counts how many copies of $H$ it occurs in.
For every set $Q$ of subgraphs of $K_n$, we define its weight as
\[
\W_Q:= \sum_{e \in \mathcal{E} (Q)} X_e = \sum_{H' \in Q} \sum_{e\in \mathcal{E}(H')} X_e
\]
Note that if an edge appears in two or more subgraphs $H'\in Q$ (copies of $H$), its weight is counted again every time.
We will let $C_H$ and $F_H$ denote the minimum weight of a partial cover and factor, respectively:
\begin{align*}
C_H(k,n) &:= \min \big\{\W_Q: Q\textrm{ is an $(k/n,H)$-cover} \big\}
\\
F_H(k,n) &:= \min \big\{\W_Q: Q\textrm{ is an $(k/n,H)$-factor} \big\}
\end{align*}
In other words, $C_H(k,n)$ (or $F_H(k,n)$) is the minimal weight of a partial cover (or partial factor) on $K_n$ that leaves at most $k$ vertices uncovered.
We will also (for technical purposes) sometimes need to keep track of upper bounds on the most expensive edge a (partial) cover or factor uses. We therefore define
\begin{align*}
C^\eps_H(k,n) &:= \min \big\{\W_Q: Q\textrm{ is an $(k/n,H)$-cover and } \max_{e\in \mathcal{E}(Q)}X_e\leq \eps \big\}\textrm{, and}
\\
F^\eps_H(k,n) &:= \min \big\{\W_Q: Q\textrm{ is an $(k/n,H)$-factor and } \max_{e\in \mathcal{E}(Q)}X_e\leq \eps  \big\}.
\end{align*}
Said cover (respectively, factor) might not exist, in which case we set $C^\eps_H(k,n)= \infty$ (resp. $F^\eps_H(k,n)=\infty$). This means that $\E C^\eps_H, \E F_H^\eps$ are not well defined, and we will simply avoid them. As we will show in the following subsection, the results from \cite{JohKahVu08}  allow us to determine a range of values of $\eps$ such that $F^\eps_H(k,n) <\infty$ (and thus $C^\eps_H(k,n) <\infty$) with very high probability. 
Note also that $C_H^\eps,F_H^\eps$ are non-increasing (random) functions of $\eps$: as $\eps$ increases, fewer edges become `forbidden', which can only decrease the minimal cost.

As every $H$-factor is also a valid $H$-cover, by definition $ C_H(k,n)\leq F_H(k,n)$ and $ C_H^\eps(k,n)\leq F_H^\eps(k,n)$.
We will also let $C_H(n):= C_H(0,n)$ and $F_H(n):= F_H(0,n)$ denote the minimal weight of an $H$-cover and $H$-factor, respectively, and similarly for $C_H^\eps(n)$, $F_H^\eps(n)$.

\section{Results and conjectures}
\label{section:thms}
Our main results are the following two theorems, where we establish bounds on $F_H$ and $C_H$, as well as prove that $F_H$ is a sharply concentrated random variable. 
\begin{theorem}
\label{mainthm}
\label{mainthm:factor}
For any graph $H$ with $d^*>1$ and $\alpha\in [0,1)$,
there are constants $0<a<b$ such that 
\[
an^{1-1/d^*} \leq F_H(\alpha n,n)\leq bn^{1-1/d^*} 
\]
with probability $1-n^{-\omega(1)}$. Furthermore, $F_H(\alpha n,n)$ is sharply concentrated around its median value $M$:
\[
|F_H(\alpha n,n)-M|= \Op(M^{3/4}).
\]
\end{theorem}
\Cref{thm:intro-main} is a special case of this theorem, with $\alpha=0$.
The two parts of the theorem are more precise versions of the statements $F_H(\alpha n,n)= \Tp(n^{1-1/d^*})$ and $F_H(\alpha n,n)/\E[F_H(\alpha n,n)]\toprob 1$, respectively.
Note, however, that together they do not guarantee that the limit (in probability) of $F_H(\alpha n,n)/n^{1-1/d^*}$ exists.
\begin{conj}
\label{conjecture}
For any graph $H$ there is a continuous decreasing function $f_H:[0,1]\to \R$ such that 
\[F_H(\alpha n,n)/n^{1-1/d^*}\toprob f_H(\alpha).\]
See \cref{guess:limit} for a discussion on what this function $f_H$ might be.
\end{conj}

 For covers, since $C_H\leq F_H$ we automatically get an upper bound by the theorem above. We also have the following lower bound.
\begin{theorem}
\label{mainthm:cover}
Let $\Delta:=\max_{H'\subset H} (\frac{e_{H'}}{v_{H'}})$. Then for any $\alpha\in[0,1)$, we have that $C_H(\alpha n,n)={\Omp (n^{1-1/\max\{d_H,\Delta\}} )}$.
\end{theorem}
Note the different exponents in the upper and lower bounds on $C_H$. They match if (for instance) $H$ is \emph{balanced}, so that $d_H=d^*$. In \cref{section:unconcentrated-cover} we discuss examples where $H$ is not balanced, only one of these bounds is sharp, and where $C_H$ is not sharply concentrated. We might still conjecture that the $H$-cover equivalent of \cref{mainthm:factor} or \cref{conjecture} holds for balanced $H$.

Before we move on to the proofs, let's briefly compare the minimum weight $H$-factor problem with the corresponding threshold problem.
In a 2008 paper, Johansson, Kahn \& Vu~\cite{JohKahVu08} found the threshold function for the appearance of an $H$-factor for strictly balanced $H$, as well as slightly less precise bounds on the threshold for general $H$. 
\begin{theorem}[Theorems 2.1 \& 2.2 in \cite{JohKahVu08}]
Assume $H$ is a fixed graph. 
\begin{enumerate}[(i)]
\item If $H$ is strictly balanced the threshold for the appearance of a $H$-factor in $G_{n,p}$ is 
$th_H:=n^{-1/d_H} (\log n)^{1/e_H}$. That is,
\[
\Pr(G_{n,p} \textrm{ contains an $H$-factor} ) =
\begin{cases}
n^{-\omega(1)}, &\textrm{ if } p\ll th_H
\\
1-n^{-\omega(1)}, &\textrm{ if } p\gg th_H
\end{cases}
\]
\item For general $H$ the threshold is $n^{-1/d^*+o(1)}$. More precisely, for any $\eps>0$,
\[
\Pr(G_{n,p} \textrm{ contains an $H$-factor} ) =
\begin{cases}
n^{-\omega(1)}, &\textrm{ if } p\ll n^{-1/d^*}
\\
1-n^{-\omega(1)}, &\textrm{ if } p\gg n^{ -1/d^*+\eps}
\end{cases}
\]
\end{enumerate} 
\end{theorem}
 This immediately implies the following upper bound on $F_H$, only a factor $n^{\eps}$ worse than the bound in \cref{mainthm:factor}. 
\begin{cor}
\label{jkv-cor}
For any $\eps>0$ and $A\gg n^{-1/d^*+\eps} $, 
$F^A_H(n) \leq  n^{1-1/d^*+\eps}$ with probability $1-n^{-\omega(1)}$.
\end{cor}

\section{Proofs}
\label{section:proofs}
In this section we state and prove several propositions from which our main theorems follow: \cref{mainthm:factor} follows from \cref{lowerbound:factor,upperbound-Hfactor,concentration}, and \cref{mainthm:cover} follows from \cref{lowerbound:cover,upperbound-Hfactor}.

\subsection{Lower bound: H-factors}
In this section we establish a lower bound on the minimum cost of $H$-factors, and then in \cref{section:lowerbound-cover} we do the same for $H$-covers. Although any lower bound on $C_H$-covers is also a lower bound on $F_H$, our lower bound for $H$-factors holds with probability $1-2^{-\Omega(n)}$, while the lower bound for $H$-covers is only shown to hold with probability $1-\eps$. For this reason we consider it worthwhile to include both.
\label{section:lowerbound}
\begin{proposition}
Assume $\alpha\in [0,1)$ is fixed (not depending on $n$).
\label{lowerbound:factor}
There exists a $c>0$ such that the minimal cost of an $(\alpha,H)$-factor is $F_H(\alpha n,n)\geq cn^{1-1/d^*}$, with probability $1-2^{-\Omega(n)}$.
\end{proposition}
To prove this, we need the following simple bound (which will also be useful several times more throughout the paper).
\begin{lemma}
\label{bound:sum-of-exps}
If $x>0$,  $X_1,X_2,\ldots X_k$ are i.i.d. $\Exp(1)$-distributed random variables and ${X:=\sum_i X_i}$, then
\[
1-x \leq \frac{\Pr(X\leq x)}{x^k/k!} \leq 1
\]
\end{lemma}
\begin{proof}
$X$ follows a Gamma distribution with shape parameter $k$ and scale parameter $1$, with density function $t^{k-1}\e^{-t}/(k-1)!$. Since $\e^{-t}\geq 1-x$ on the interval $t\in [0,x]$, 
\begin{equation}
    \Pr(X\leq x)\geq (1-x)\int_0^x \frac{t^{k-1}}{(k-1)!}dt = (1-x)x^k/k!.
\end{equation}
Similarly, using $\e^{-t}\leq 1$ gives $\Pr(X\leq x)\leq x^k/k!$.
\end{proof}
\medskip

 We can now prove \cref{lowerbound:factor}.

\begin{proof}[Proof of \cref{lowerbound:factor}]
Assume without loss of generality that $\alpha n$ is an integer multiple of $v_H$. Let $t$ be the smallest number of copies of $H$ an $(\alpha,n)$-factor can have. Since $(1-\alpha) n$ vertices of $K_n$ are covered, each by a unique copy of $H$, $v_Ht= (1-\alpha) n$.

We will first prove that $F_H(\alpha n,n)\geq cn^{1-1/d_H}$ whp by applying a first moment method to the following random variable.
For any $L=L(n)$, let $Y_L$ be the number of $(\alpha,H)$-factors $Q$ that have precisely $t$ copies of $H$ and that have a weight $W_Q\leq L$. Note that if $Y_L=0$ then $F_H(\alpha n,n)>L$, because any $(\alpha,H)$-factor that has more than $t$ copies of $H$ contains one with precisely $t$ copies.

How many $(\alpha,H)$-factors in $K_n$ with precisely $t$ copies of $H$ are there (regardless of weight)?
There are $\binom{n}{\alpha n}= 2^{O(n)}$ ways to pick which $\alpha n$ vertices will not be covered, and then at most $(v_H t)!/t!=2^{O(n)}n^{(v_H-1)t}$ ways to construct an $H$-factor on the remaining $v_H t=(1-\alpha) n $ vertices. Rewriting the exponent of $n$ as $v_H-1=e_H/d_H$, we can upper bound the number of such factors by $(c_1 n^{1/d_H})^{e_H t}$ for some constant $c_1$.
Consider now an $(\alpha,H)$-factor $Q$ with $t$ copies of $H$. It consists of $e_H t $ edges, so by \cref{bound:sum-of-exps} 
\begin{equation}
\label{c2}
    \Pr(W_Q\leq L)\leq L^{e_H t}/(e_H t)! \leq   (c_2L/n)^{e_H t},
\end{equation}
for some $c_2>0$. We therefore get that
\(
\E Y_L\leq  (c_1 c_2 L n^{-1+1/d_H})^{e_Ht}
\).
Since $c_1,c_2$ are constants, we can ensure that the expression within brackets is at most $1/2$ by letting $L:=c n^{1-1/d_H}$ for a sufficiently small $c=c(\alpha, H)$.
Then $\E Y_L\leq  2^{-e_H t}=2^{-\Omega(n)}$, whence $F_H(\alpha n,n) \geq c n^{1-1/d_H}$ with probability $2^{-\Omega(n)}$.

Now, if $d^*>d_H$ we can improve this lower bound. 
Let $H^*\subseteq H$ be a subgraph of the maximal density $d^*$. Consider $Q$ as above: an $(\alpha,H)$-factor which consists of $t$ copies of $H$, with $v_H t=(1-\alpha)n$. This partial $H$-factor will contain a partial $H^*$-factor $Q^*$ consisting of $t$ copies of $H^*$ and hence covering $tv_{H^*}$ vertices -- just remove the superfluous vertices and edges from each copy of $H$ in $Q$. 
This $Q^*$ is an $(\alpha^*,H^*)$-factor, with $\alpha^*$ such that number of vertices covered by $Q^*$ is $(1-\alpha^*)n=v_{H^*}t=\Omega(n)$. 
By the previous argument (and since $\alpha^*\in [0,1)$) $F_{H}(\alpha,n)\geq F_{H^*}(\alpha^*,n) \geq c(\alpha^*,H^*) n^{1-1/d^*}$ with probability $2^{-\Omega(n)}$. 
\end{proof}

In the following remark we discuss some possible optimizations of this result.
\begin{remark}
\label{guess:limit}
With some more care taken, we can find minimal $c_1,c_2$ in the proof above. The number of $H$-factors is $n!/(\alpha n)! t!Aut(H)^t$ (where $Aut(H)$ is the number of automorphisms of $H$). Applying Stirling's approximation to this and to $(e_Ht)!$ in (\ref{c2}) leads to $c_1c_2=\frac{r}{e_H} \e^{1-1/d_H} \cdot(r \alpha^{-\alpha r}Aut(H))^{-1/e_H} $, where ${r:=n/t}={v_H/(1-\alpha)}$. 
It is a tempting conjecture that the resulting bound with $c^{-1}:= c_1c_2$ is tight, at least for strictly balanced $H$. In other words, that $F_H(n)/n^{1-1/d_H}$ should converge in probability to this $c$.
\end{remark}
\subsection{Lower bound: H-covers}
\label{section:lowerbound-cover}
 We now prove the less sharp lower bound on the minimal cost of an $H$-cover.
\begin{proposition}    
\label{lowerbound:cover}
For any fixed $\alpha>0$ there exists a $K>0$ such that for any $t>0$ fixed or tending to $0$ as $n\to \infty$,
\begin{enumerate}[(i)]
    \item $C_H(\alpha, n)\geq t n^{1-1/d_H}$ with probability at least $1-Kt^{e_H}$.
    \item Let $\Delta:= \max_{G\subseteq H}(e_G/v_G)$. Then $C_H(\alpha, n)\geq t n^{1-1/\Delta}$ with probability at least $1-Kt^{e_G}$, where $G$ is the graph that attains the maximum $\Delta$.
\end{enumerate}
\end{proposition}
\begin{proof}
For any $b>0$, call a copy $H'\subset K_n$ of $H$ \emph{$b$-cheap} if $\W_{\{H'\}}< b$, i.e. if the total weight of the edges in $H'$ is at most $b$.
Let $N_b$ be the total number of $b$-cheap $H'$. We want to estimate $\E[N_b]$.
For a given $H'$, by \cref{bound:sum-of-exps} the probability that it is $b$-cheap is at most $b^{e_H}/{e_H}!$.
Furthermore, there are less than $n^{v_H}$ copies of $H$ in $K_n$.
Then, by Markov's inequality, for any $\lambda>0$,
\begin{equation}
    \Pr\left(N_b\geq \lambda\right) \leq \frac{\E[N_b]}{\lambda} \leq \frac{n^{v_H}b^{e_H}}{\lambda{e_H}!}
    \label{number-of-cheap-H}
\end{equation}
Now, suppose that there exists an $(\alpha,H)$-cover $Q$ with $\W_Q\leq tn^{1-1/d_H}$.
This $Q$ consists of at least $\frac{\alpha}{v_H}n$ copies of $H$, since each copy of $H$ covers at most $v_H$ vertices not covered by another copy.

The number of $H'\in Q$ that are not $b$-cheap can be at most $\W_Q/b$. In particular for $b:={2v_H}n^{-1/d_H}/{\alpha}$, there can be at most ${\alpha n}/{2v_H}$ that are not $b$-cheap, or in other words at most half of the $H'\in Q$. Hence $Q$ must contain at least $\frac{\alpha}{2v_H}n$ $b$-cheap copies $H'$, which implies that $N_b\geq \frac{\alpha}{2v_H}n$.
By (\ref{number-of-cheap-H}), 
\begin{equation}
   \Pr(N_b\geq \frac{\alpha}{2v_H}n)\leq \frac{2v n^{v_H}b^{e_H}}{\alpha n{e_H}!}=\frac{(2vt/\alpha)^{e_H}}{\alpha {e_H}!} 
\end{equation}
This immediately implies part (i).
For part (ii), 
consider the subgraph $G$ that attains the maximum $\Delta:= \max_{G\subset H}(e_G/v_G)$.
As noted earlier, any $(\alpha,H)$-cover $Q$ contains at least $\frac{\alpha}{v_H}n$ copies of $H$. Let $H_1,H_2, \ldots$ be an enumeration of them, and $G_i\subset H_i$ be copies of $G$ in each. Note that we might have $G_i=G_j$ for some $i\neq j$, as two distinct copies of $H$ might overlap in a copy of $G$. 
\begin{equation}
\label{ineq:WQ}
   \W_Q=\sum_i \W_{H_i} \geq \sum_i \W_{G_i} \geq \frac{\alpha}{v_H}n \min \W_{G'}, 
\end{equation}
where the last minimum is taken over all copies $G'\subset K_n$ of $G$.
Applying (\ref{number-of-cheap-H}) with $\lambda=1$, $G$ instead of $H$ and some $b$ to be determined shortly, we see that 
$\Pr\left(N_b\geq 1\right)$
is at most 
${(nb^{\Delta})^{v_G}}/{e_G!}$.
Letting $b=t n^{-1/\Delta}$ for a small $t>0$, the right-hand side of (\ref{ineq:WQ}) is at most $t^{v_G}/e_G!$.
Hence $\W_Q\geq \alpha n\min \W_{G'}/v_H\geq t \alpha n^{1-1/\Delta}/v_H$ with probability at least  $1-t^{v_G}/e_G!$, from which (ii) follows.
\end{proof}
\begin{remark}
For strictly balanced $H$, \cref{lowerbound:cover}(i) can be sharpened by a second moment argument to hold with probability $1-o(1)$ rather than $1-Kt^{e_H}$. 
\end{remark}

\subsection{Red-green split lemma}
In this section we introduce the red-green split trick mentioned in \cref{subsection:completingpartialfactors}.
This lemma will be useful both to prove the upper bound on $F_H$, as well as to prove that it is sharply concentrated. It is also used in \cref{pseudo-dim}.

We state and prove \cref{redgreensplit} (as well as \cref{upperbound-Hfactor}) not only for $F_H$ but for $F_H^A$: the minimum weight of an $H$-factor using no edge of weight more than $A$, considering then $F_H$ as the particular case where $A=\infty$. Keeping track of upper bounds on the most expensive edge in an $H$-factor make statements and proofs slightly more involved.
While such bounds will be of use in \cref{thm:pseudo-dim}, they are not necessary for our main results, \cref{mainthm:cover,mainthm:factor}.
We therefore suggest that the reader who is only interested in the latter theorems simply ignore the superscript in $F_H^A$, and any inequalities involving $A, A_k, B $ and $C$.

\label{section:redgreensplit}
\begin{lemma}
\label{redgreensplit}
Let $n>m>k\geq 0$ be integer multiples of $v_H$.
\begin{enumerate}[(1)]
    \item For any $t\in (0,1)$, the random variables $F^A_H(m,n)$, $F_H^B(k,m)$ and $F_H^C(k,n)$ (where $C\geq \max(\frac{A}{t},\frac{B}{1-t})$) can be coupled such that surely
    \[
    F_H^C(k,n)\leq \frac{F^A_H(m,n)}{t}+\frac{F_H^B(k,m)}{1-t}.
    \]
    \item Let $a,b,A,B>0$ and let $C\geq (a+b)\max(A/a,B/b)$. Then
    \[
    \Pr\Big(F^C_H(k,n)> (a+b)^2\Big)
    \leq \Pr\Big(F^A_H(m,n)> a^2\Big)
    + \Pr\Big(F^B_H(k,m)> b^2\Big).
    \]
\end{enumerate}
Both of these inequalities also hold when $A=B=C=\infty$, i.e.\ with $F_H$ instead of $F_H^A$, $F_H^B$ and $F_H^C$.
\end{lemma}
\begin{remark}
The lemma also holds for $H$-covers, and in that case the requirement that $n,m$ and $k$ are integer multiples of $v_H$ is not necessary. The proof for $H$-covers is mutatis mutandis. However, we will only prove and use the lemma for factors.  
\end{remark}

\begin{proof}
We will begin by proving part (1) of the lemma.
Let $G$ be the multigraph on $[n]$ given by connecting every pair of vertices by two parallel edges, one green and one red. Independently for all edges, assign each green edge an $\Exp(t)$-distributed random weight and each red edge an $\Exp(1-t)$-distributed random weight.
We will use the following properties of the exponential distribution:
\begin{enumerate}[(i)]
    \item if $X\sim \Exp(t)$ and $Y\sim \Exp(1-t)$ are independent, then $\min(X,Y)\sim \Exp(1)$
    \item if $X\sim \Exp(t)$, then $tX\sim \Exp(1)$.
\end{enumerate}

 Let $Z$ be the cost of the cheapest $(k/n,H)$-factor in $G$ that uses no edge more expensive than $C$. (If no such factor exists, $Z=\infty$.) It will always use the cheaper one of two parallel edges, so by property (i) we see that $Z\overset{d}{=} F^C_H(k,n)$.
Our aim is now to construct a fairly cheap (but not necessarily optimal) such factor in $G$. First, we pick the cheapest green $(m/n,H)$-factor that uses no edge more expensive than $A/t$, and let $Z_{green}$ be its cost. Note that by the rescaling property (ii), $tZ_{green}\overset{d}{=} F^A_H(m,n)$.

We are left with a random set of $m$ uncovered vertices. Crucially, this random set is independent from the weights on the red edges. Pick the cheapest red $(k/m,H)$-factor (i.e. a partial factor leaving at most $k$ out of $m$ vertices uncovered) on this set that uses no edge more expensive than $B/(1-t)$, and let its cost be $Z_{red}$. Again by (ii), $ (1-t)Z_{red}\overset{d}{=} F^B_H(k,m)$.

Combining the green copies of $H$ from the first step with the red copies of $H$ in the second step gives us a partial $H$-factor $Q$ on $G$ covering all but at most $k$ vertices -- i.e.\ a $(k/n,H)$-factor.
No edge in $Q$ costs more than $\max(\frac{A}{t},\frac{B}{1-t})\leq C$, whence $Z\leq W_Q = Z_{green}+Z_{red}$.
Thus (by an appropriate coupling) the following inequality holds:
\begin{equation}
    F^C_H(k,n)\leq \frac{F^A_H(m,n)}{t}+\frac{F^B_H(k,m)}{1-t}.
\end{equation}
For part (2), it follows from part (1) that if $F^A_H(m,n)\leq a^2$ and $F^B_H(k,m)\leq b^2$, then $F^C_H(k,n) \leq \frac{a^2}{t}+\frac{b^2}{1-t}$. Minimizing over $t$ gives that the right hand side is $(a+b)^2$ for $t=a/(a+b)$, and for this $t$ we get that $C=\max(\frac{A}{t},\frac{B}{1-t})=(a+b)\max(A/a,B/b)$. 
Hence $F^C_H(k,n)\leq (a+b)^2$, unless $F^A_H(m,n)> a^2$ or $F^B_H(k,m)> b^2$. Using the union bound on these two events give the inequality in part (2).

For the case $A=B=C=\infty$ the proof is nearly identical, except we do not need to keep track of the cost of the most expensive edges.
\end{proof}

\subsection{Upper bound}
In this section, we prove the following upper bound on the total cost of an $H$-factor, both unconstrained and limited to using only edges of weight at most $A$.

\label{section:upperbound}
\begin{proposition}
\label{upperbound-Hfactor}
For any fixed graph $H$ with $d^*>1$ and any $\eps>0$, there exists a $c>0$ such that if $A\geq n^{-1/d^*+\eps}$, then  $F^A_H(n)\leq cn^{1-1/d^*}$ with probability at least ${1-n^{-\omega(1)}}$. In particular, this holds for $A=\infty$.
\end{proposition}
To prove this proposition, we will need the following theorem from \cite[Thm 4.9]{jsl-randomgraphs}, originally due to Rucinski~\cite{rucinski-partialfactor}.
\begin{theorem}
\label{upperbound-partialfactor}
\label{gnp-has-partial-factor}
For any $\alpha\in (0,1)$ there exist constants $c,t>0$ such that $G_{n,p}$ with $p=cn^{-1/d^*}$ contains an $(\alpha,H)$-factor with probability at least $1-2^{-tn}$.
\end{theorem}
In \cite{jsl-randomgraphs}, the existence of such a partial factor is only stated to hold with probability $1-o(1)$, but in the proof the probability is shown to be $1-2^{-\Omega(n)}$.

\begin{proof}[Proof of \cref{upperbound-Hfactor}]
The proof strategy is essentially this: For some small fixed number $\alpha>0$, we will find a cheap $H$-factor on $n$ vertices by iteratively using the red-green split trick from \cref{redgreensplit}. This will give a cheap $(\alpha,H)$-factor on $n_i$ vertices (starting with $n_0:=n$), then a cheap $(\alpha,H)$-factor on the remaining $n_{i+1}$ vertices, and so on, for a total of $k$ steps. On the remaining $n_k$ vertices, it suffices to find a not too expensive $H$-factor.

More precisely, pick $\alpha$ so that $\alpha^{1-1/d^*}= \frac{1}{4}$ (and hence $\alpha< \frac{1}{4}$).
Let $n_0:=n$ and let $n_{i}$ be the largest multiple of $v_H$ such that $n_{i}\leq \alpha n_{i-1}$.
Also for some small fixed $\delta>0$ to be determined later, let $k$ be an integer such that $\alpha^k \leq n^{-\delta} \leq \alpha^{k-1}$.
For this choice of $n_i$ and $k$, we have that $\alpha^{i+1}n \leq n_i \leq \alpha^i n$ and $\alpha n^{1-\delta}\leq n_k \leq  n^{1-\delta}$. Also, $4^k<n^{\delta}$. 

Applying part (1) of \cref{redgreensplit}, with $t=1/2$ and $A_i:=2^i A$, repeatedly to $F^{A_{i}}_H(n_{i})$ for $i=0,1,\ldots ,{k-1}$, we  get that there exists a coupling such that 
\begin{align}
\label{redgreenrecursion}
F_H^{A_0}(n_0) &\leq 2F^{A_{1}}_H(n_1,n_0)+2F^{A_{1}}_H(n_{1}) \nonumber
\\
&\leq 2F^{A_{1}}_H(n_1,n_0)+4F^{A_{2}}_H(n_2,n_1)+4F^{A_{2}}_H(n_{2}) \nonumber
\\
&\leq 2F^{A_{1}}_H(n_1,n_0)+4F^{A_{2}}_H(n_2,n_1)+8F^{A_{3}}_H(n_3,n_2)+8F^{A_{3}}_H(n_{3}) \nonumber
\\
&\ldots \nonumber
\\
& \leq \underbrace{\sum_{i=0}^{k-1} 2^{i+1} F^{A_{i+1}}_H\Big(n_{i+1},n_{i}\Big)}_{(\ref{redgreenrecursion}a)} +  \underbrace{2^{k}F^{A_k}_H(n_k)}_{(\ref{redgreenrecursion}b)}
\end{align}
Let's begin with the sum (\ref{redgreenrecursion}a).
By \cref{upperbound-partialfactor}, there exists constants $c,t$ (depending only on $\alpha, H$) such that if $A_{i+1}\geq  cn_i^{-1/d^*}$ then $F^{A_{i+1}}_H(n_{i+1},n_i)\leq cn_i^{1-1/d^*}$ with probability at least $1- 2^{-tn_i}\geq 1- 2^{-tn_k}$.
To check whether this lower bound on  $A_{i+1}$ holds, note that since $A_i\geq A$ and $n_i\geq n_k$, it suffices to show that $A n_k^{1/d^*}\geq  C$.
Using that $n_k \geq \alpha^{k+1} n$ and $\alpha^k\geq \alpha n^{-\delta}$, we get that
\begin{equation}
A n_k^{1/d^*}
= n^{-1/d^*+\eps}n_k^{1/d^*}
\geq   n^{\eps} (\alpha^{k+1})^{1/d^*}
\geq  n^{\eps}  (\alpha^2 n^{-\delta} )^{1/d^*}
\gg n^{\eps/2} ,
\label{bound:Ank}
\end{equation}
where the last inequality holds by picking $\delta$ sufficiently small.
Hence the conditions of \cref{upperbound-partialfactor} are met, and it then follows (by a union bound) that with probability at least $1-k 2^{-tn_k}= 1-n^{-\omega(1)}$, we have (\ref{redgreenrecursion}a) $\leq  2c\sum_{i=0}^{k-1} 2^i n_i^{1-1/d^*}$.
Since $n_i\leq \alpha^i n$ and $\alpha^{1-1/d^*}= \frac{1}{4}$ (by the choice of $\alpha$), we can bound the terms in this sum by
\begin{equation}
\label{bound:2ini}
2^i n_i^{1-1/d^*}\leq (2\alpha^{1-1/d^*})^i \cdot n^{1-1/d^*}\leq 2^{-i} n^{1-1/d^*}.
\end{equation}
Hence (\ref{redgreenrecursion}a) is at most $ 4cn^{1-1/d^*}$ whp.
For the term (\ref{redgreenrecursion}b) of equation (\ref{redgreenrecursion}), the slightly rougher bound in \Cref{jkv-cor} suffices: for any $\delta'>0$, if $A_k\gg n_k^{1-1/d^*+\delta'}$ then $F^{A_k}_H(n_k)\leq n_k^{1-1/d^*+\delta'}$ with probability $n_k^{-\omega(1)}$.
But by (\ref{bound:Ank}), $A_k \geq A \gg n_k^{-1/d^*+\eps/2}$, so the condition on $A_k$ is met if we pick $\delta'<\eps/2$.
Then
\begin{equation}
(\ref{redgreenrecursion}b) = 2^k F^{A_k}_H(n_k)
\leq 2^k  n_k^{1-1/d^*+\delta'}
\leq 2^{-k} n^{1-1/d^*+\delta'}
\end{equation}
where the last inequality uses inequality (\ref{bound:2ini}) and $n_k^{\delta'}\leq n^{\delta'}$.
From the choice of $k$, $2^{-k}\leq n^{-\delta /|\log_2\alpha|}$, and we can therefore ensure that the right-hand side above is $o(n^{1-1/d^*})$ by picking $\delta'$ sufficiently small ($\delta'<\delta /|\log_2\alpha|$).
It follows that $(\ref{redgreenrecursion}a) +(\ref{redgreenrecursion}b) \leq {(4c+o(1))n^{1-1/d^*}}$ with probability $1-n^{-\omega(1)}$.
\end{proof}

\subsection{Concentration}
\label{section:concentration}
We will now move on to show that $F_H$ is sharply concentrated.
\begin{proposition}
\label{concentration}
For any graph $H$ with $d_H>1$, $\eps >0$ and $\alpha \in [0,1)$, there exists a $c>0$ such that if we let $M=M(\alpha,n,H)$ denote the median of $F_H(\alpha n,n)$, then for all sufficiently large $n$ and with probability at least $1-\eps$,
\[
|F_H(\alpha n,n)-M|<cM^{3/4}.
\]
\end{proposition}
Let's consider a dual problem: How large is the largest partial factor that costs at most $L$, for some some $L=L(n)$?
More precisely, let the random variable $Z_H=Z_H(n,L)$ be defined by
\[Z_H:=\max\{\alpha n: \textrm{there exists an } (1-\alpha,H)\textrm{-factor } Q \textrm{ with }\W_Q\leq L\}.\]
In other words, $Z_H$ is the largest number of vertices that a partial factor costing at most $L$ can cover. Note that $Z_H(n,L)\geq n-m$ if and only if $F_H(m,n) \leq L$.
Our first step is to apply Talagrand's concentration inequality to $Z_H$. To do so we need the definitions of $f$-certifiable and Lipschitz random variables.
\begin{definition}[$f$-certifiable random variable]
Let $X:\Omega^n \to \R$ be a random variable. For a function $f$ on $\R$ we say that $X$ is $f$-certifiable if for any $\omega\in \Omega^n$ with $X(\omega)\geq s$, there is a set $I\subseteq [n]$ of at most $f(s)$ coordinates such that $X(\omega')\geq s$ for all $\omega'$ which agree with $\omega$ on $I$.
(That is, $\omega'_i=\omega_i$ for all  $i \in I$.)
\end{definition}

\begin{definition}[Lipschitz random variable]
Let $X$ be as above. We say that $X$ is $K$-Lipschitz if for every $\omega,\omega'$  with $\omega_i=\omega_i'$ for all but one $i$, $|X(\omega)-X(\omega')|\leq K$. 
\end{definition}
We can now state Talagrand's inequality. While it was first established in~\cite{talagrand1995concentration}, we will use the following, more `user-friendly', version from~\cite{probmeth}.
\begin{theorem}[Talagrand's concentration inequality]
\label{talagrand}
Assume $\Omega$ is a probability space. If $X$ is a $K$-Lipschitz, $f$-certifiable random variable $X:\Omega^n\to \R$ where $\Omega^n$ is equipped with the product measure, then for any $b,t\geq 0$,
\[
\Pr(X\leq b)\Pr(X\geq b+tK\sqrt{f(b)})\leq \exp({-t^2/4}).
\]
\end{theorem}
 The following lemma finds the appropriate values of $f$ and $K$ so that we can apply this inequality to the random variable $Z_H$.
\begin{lemma}
$Z_H$ is $v_H$-Lipschitz and $f$-certifiable with $f(s)=e_H\lceil\frac{s}{v_H} \rceil\leq \frac{e_H}{v_H} n$.
\end{lemma}
\begin{proof}[Proof]
To show that $Z_H$ is $f$-certifiable, pick an integer $s\in [n]$ and a tuple of edge weights  $\omega\in \Omega^{\binom{n}{2}}$ such that $Z_H(\omega)\geq s$.
Then there exists a partial $H$-factor $Q$ with $\W_Q(\omega)\leq L$ and which covers at least $s$ vertices. Assume without loss of generality that $Q$ is one of the smallest such partial $H$-factors.
It then contains $\lceil s/v_H\rceil$ copies of $H$ and $f(s):=e_H \lceil s/v_H\rceil$ edges. For any $\omega'$ which agree with $\omega$ on the $f(s)$ edges of $Q$, $\W_Q(\omega')=\W_Q(\omega)\leq L$. Hence $Z_H(\omega')\geq s$. (It might be that $Z_H(\omega)\neq Z_H(\omega')$, here we only care whether they are $\geq s$.)

To show the Lipschitz condition, pick an edge $e$ and condition on all other edge weights. Consider $Z_H$ as a function of just $x=X_e$.
Note first that $Z_H(x)$ is a non-increasing function, i.e. $Z_H(x)\leq Z_H(0)$ for any $x\geq 0$.
Let $Q$ be a partial $H$-factor achieving the maximum size $Z_H(0)$. That is, $Q$ covers $Z_H(0)$ vertices and has weight $W_Q=W_Q(x)$ such that $W_Q(0)\leq L$.
Is $e\in \mathcal{E}(Q)$?
\begin{enumerate}[(i)]
\item If $e\in \mathcal{E}(Q)$, let $H_e$ be the copy of $H$ in $Q$ which contains $e$. Then $Q-H_e$ is a partial $H$-factor with weight at most $W_{Q-H_e}(x)< W_Q(0)\leq L$ (for any  $x$), and it covers $Z_H(0)-v_H$ vertices. Hence $Z_H(x) \geq Z_H(0)-v_H$.
\item If $e\notin \mathcal{E}(Q)$, then $W_Q(x)$ is a constant function and $W_Q(x)=W_Q(0)\leq L$. Hence $Z_H(x) \geq Z_H(0)$.
\end{enumerate}
In either case, $Z_H(0)-v_H\leq Z_H(x)\leq Z_H(0)$. Thus $Z_H$ is $v_H$-Lipschitz.
\end{proof}

\begin{remark}
This is where our proof would fail for the corresponding cover problem. For covers, some edge might belong to a large number of copies of $H$, leading to a large Lipschitz constant. This is the case in our example in \cref{section:unconcentrated-cover}.
\end{remark}

Before proceeding with the proof of \cref{concentration}, we'll need two small lemmas.
\begin{lemma}
\label{monotone-cost}
If $k<m<n$, $F_H(m,n)\leq \frac{n-m}{n-k} F_H(k,n)$.
\end{lemma}
 \begin{proof}
$F_H(m,n)$ is the lowest cost of a partial $H$-factor covering at least $n-m$ vertices of $K_{n}$. We can construct a cheap such partial factor in two steps: First, let $Q$ be the optimal $(\frac{k}{n},H)$-factor (which consists of $(n-k)/v_H$ copies of $H$), i.e. the factor such that $\W_Q=F_H(k,n)$.

Next, let $Q'$ be the partial factor obtained by removing all but the ${(n-m)/v_H}$ cheapest copies of $H$ in $Q$, leaving a $(\frac{m}{n},H)$-factor. Then $Q'$ contains a fraction ${(n-m)/(n-k)}$ of the copies of $H$ in $Q$. Hence it costs ${\W_{Q'}\leq \frac{n-m}{n-k}\W_Q}$.
\end{proof}

\begin{lemma}
\label{FH-is-continuous}
For any $m$ and $n$, the random variable $F_H(m,n)$ follows a continuous distribution (i.e. it has no atoms).
\end{lemma}
 \begin{proof}
For any partial factor $Q$ and $t\geq 0$, $\Pr(\W_Q=t)=0$. And since there are only finitely many such $Q$, $\Pr(F_H(m,n)=t)\leq \Pr(\exists Q: \W_Q=t)=0$.
\end{proof}

We can now finally prove that the cost of a (partial) $H$-factor concentrates around its median.

 \begin{proof}[Proof of \cref{concentration}]
Let $m$ be the largest multiple of $v_H$ such that $m\leq  \alpha n $.
By \cref{FH-is-continuous}, $F_H(m,n)$ is a continuous random variable, whence we can find $L$ such that $\Pr(F_H(m,n)\leq L) = \eps$. Using the upper bound (\cref{upperbound-Hfactor}) and lower bound (\cref{lowerbound:factor}) on $F_H$, we see that in order for $\Pr(F_H(m,n)\leq L) = \eps$ to hold, we must have $L=\Theta(n^{1-1/d^*})$. (For the lower bound, the condition $\alpha <1$ is used.)
Now, let's apply the Talagrand inequality to the $v_H$-Lipschitz, ${e_H n/v_H}$-certifiable random variable $Z_H(L,n)$. Choose $t>0$ such that $\exp({-t^2/4}) = \eps^2$ and let $b:= n-m-k$, where $k:=\lceil t\sqrt{e_Hv_H n}\rceil$. Then 
\begin{equation}
\label{ZH-talagrand}
\Pr(Z_H\leq n-m-k)\cdot \Pr(Z_H\geq n-m) \leq \eps^2.
\end{equation}
By the choice of $L$ and recalling that $Z_H(L,n)$ is the largest $n-m$ such that $F_H(m,n)\leq L$, the second probability in the left-hand side of (\ref{ZH-talagrand}) is $\eps$. Hence the first probability is
\begin{align}
\label{FH-talagrand}
\Pr(F_H(m+k,n)\geq L)&=\Pr(Z_H\leq n-m-k)\leq \eps.
\end{align}
So with probability at least $1-\eps$, there is a partial $H$-factor of cost at most $L$ and that leaves at most $m+k$ vertices uncovered. 
What is the cost of a partial factor covering $k$ out of the remaining $m+k$ vertices? By \cref{monotone-cost} and \cref{upperbound-Hfactor}
\begin{equation}
\label{ineq:smallpatch}
    F_H(m,m+k)
\leq \frac{k}{m+k}F_H(m+k)
\leq  ck (m+k)^{-1/d^*}
\leq  ck^{1-1/d^*}=:\ell,
\end{equation}
with probability $1-k^{-\omega(1)}\geq 1-\eps$ for some constant $c$ (since $k\gg 1$).
Using part (2) of \cref{redgreensplit},
\begin{align}
\label{bound:FH-L}
\Pr\left(F_H(m,n) >\big(\sqrt{L}+\sqrt{\ell}\big)^2\right)
\leq \phantom{+} &\Pr(F_H(m+k,n) > L) 
\\
\label{bound:FH-l}
+&\Pr(F_H(m,m+k) >\ell)
\end{align}
Note that $\ell=\Theta(\sqrt{L})$, since $L=\Theta(n^{1-1/d^*})$, $\ell=\Theta(k^{1-1/d^*})$ and $k=\Theta(\sqrt{n})$.
Thus 
$\big(\sqrt{L}+\sqrt{\ell}\big)^2 \leq L + bL^{3/4}$ for some constant $b$. 
The right-hand side of (\ref{bound:FH-L}) is at most $\eps$ by (\ref{FH-talagrand}), and (\ref{bound:FH-l}) is at most $\eps$ by (\ref{ineq:smallpatch}).
Thus (\ref{bound:FH-L}),(\ref{bound:FH-l}) give that 
\begin{align}
&\Pr(F_H(m,n) > L+bL^{3/4}) \leq 2\eps,
\end{align}
and by the choice of $L$, $\Pr(F_H(m,n) < L) = \eps$. Assuming without loss of generality that $\eps<1/4$, this also implies that the median $M$ of $F_H(m,n)$ lies in the interval $[L,L+bL^{3/4}]$, and in particular $M=\Theta(L)$. Hence $|F_H(m,n)-M| =\Op(M^{3/4})$.
\end{proof}

\section{Other edge weight distributions}
As mentioned in \cref{intro:pseudodim}, the exact edge weight distribution doesn't matter, only its asymptotic behavior near $0$. Here we prove this fact.
\label{pseudo-dim}
\begin{theorem}
\label{thm:pseudo-dim}
Assume $K_n$ is equipped with positive i.i.d.\ edge weights $Z_e$ with some common cdf $\tilde G$ satisfying $\lim_{x\to 0} \tilde G(x)/x=1$ (i.e.\ $\tilde G(x) = x+o(x)$). Let $\tilde F_H(m,n)$ be the minimum weight $(m/n,H)$-factor with respect to these weights (and similarly for $\tilde F_H(n)$, $\tilde C_H(n)$, $\tilde C_H(m,n)$). Then these edge weights can be coupled to i.i.d. $\Exp(1)$ edge weights in such a way that for any $m=m(n)$ with $\lim_{n\to \infty} m/n <1$, $\tilde F_H(m,n)/F_H(m,n)\toprob 1$ and $\tilde C_H(m,n)/C_H(m,n)\toprob 1$.
\end{theorem}
\begin{remark}
If instead $\tilde G(x)=\lambda x+o(x)$ for some $\lambda>0$, we can replace the edge weights $Z_e$ with weights $\lambda Z_e$. This changes the optimal cost by a factor $\lambda$, and since $\Pr(\lambda Z_e \leq x)=\tilde G(x/\lambda) = x+o(x)$, we have that $\tilde F_H(m,n)/F_H(m,n)\toprob \lambda$.
\end{remark}

 \begin{proof}[Proof of \cref{thm:pseudo-dim}]
We will prove this for $m=0$ and $F_H$ -- the proof is essentially identical for $m>0$ and/or $C_H$, but the notation becomes messier.

Let $X_e\sim \Exp(1)$, and let $G(x)=1-\e^{-x}$ be the CDF of this distribution.
Then $G(X_e)$ is uniformly distributed in the interval $[0,1]$, and we can therefore couple it to $Z_e$ by letting $Z_e:=\tilde G^{-1}(G(X_e))$.

Pick a small fixed $\eps>0$. 
Since both $\tilde G(x)$ and $G(x)$ are asymptotically $x+o(x)$ as $x\to 0$, we can find a $C=C(\eps)>\eps$ such that for any $x\in [0,3C]$, both
$G(x)\leq \tilde G((1+\eps)x)$ and $\tilde G(x)\leq  G((1+\eps)x)$ holds.
So whenever either $X_e$ or $Z_e$ is at most $2C$, the other is at most $2C(1+\eps)<3C$, and hence
\begin{equation}
(1-\eps) X_e\leq Z_e\leq (1+\eps)X_e. \label{weights-sandwich}
\end{equation}
We will prove that the following chain of inequalities hold whp:
\begin{equation}
\label{ineqs:sandwich}
1-4\eps\leq \frac{\tilde F_H(n)}{ F_H^{2C}(n)}\leq  1+\eps
\end{equation}
For the second inequality of (\ref{ineqs:sandwich}), consider $F_H^{2C}(n)$. This is finite iff there exists an $H$-factor $Q$ which uses no edge of weight more than $2C$. 
We know from \cref{jkv-cor} that such $Q$ exists with probability $1-n^{-\omega(1)}$, so it is enough to prove (\ref{ineqs:sandwich}) holds whp under the assumption that there is such $Q$, or equivalently, that $F_H^{2C}(n)< \infty$. Pick $Q$ as the cheapest such $H$-factor, so that $\W_Q=F_H^{2C}(n)$.  An edge $e$ in $Q$ has edge weight $X_e\leq 2C$ by construction, whence $Z_e\leq (1+\eps)X_e$, and
\begin{equation}
    \tilde F_H(n)\leq \sum_{e\in \mathcal{E}(Q)}Z_e \leq  \sum_{e\in \mathcal{E}(Q)}(1+\eps)X_e = (1+\eps) F_H^{2C}(n).
\end{equation}
For the first inequality of (\ref{ineqs:sandwich}), let $Q$ instead be the optimal $H$-factor with respect to the edge weights $Z_e$, i.e.\ $\sum_{e\in \mathcal{E}(Q)}Z_e=\tilde F_H(n)$. We will use it to construct a cheap $H$-factor (w.r.t.\ $X_e$).
Call a copy $H'\in Q$ `bad' if it contains at least one edge $e$ with cost $Z_e\geq C$. The total number of such edges in $Q$ is at most $\tilde F_H(n)/C$, so there are at most this many bad copies, and at most $v_H\tilde F_H(n)/C$ vertices are covered by a bad copy.

Using the second inequality of (\ref{ineqs:sandwich}) together with \cref{upperbound-Hfactor}, we see that $\tilde F_H(n)\leq (1+\eps)F^{2C}_H(n) \leq Kn^{1-1/d^*}$ whp for some constant $K$, and then at most $k:= {v_H \cdot \lfloor K n^{1-1/d^*}/C\rfloor}\ll n$ vertices are covered by a bad copy. Removing every bad copy then gives an $(\frac{k}{n},H)$-factor using no edge more expensive than $C$ (whp). Hence $\tilde F_H(n) \geq (1-\eps) F^C_H(k,n)$. By \cref{redgreensplit}, 
\begin{equation}
\label{ineq:pseudodim-redgreen}
F^{2C}_H(n)
\leq 
\frac{F^C_H(k,n)}{1-\eps} + \frac{F^{C\eps}_H(k)}{\eps}.  
\end{equation}
Pick some $a_n,b_n$ with $k^{1-1/d^*}\ll a_n\ll b_n \ll  n^{1-1/d^*}$. The second term on the right-hand side of (\ref{ineq:pseudodim-redgreen}) is by \cref{upperbound-partialfactor} at most $a_n$ whp. On the other hand, by \cref{lowerbound:cover} the first term is at least $b_n$ whp. Hence $F^{C\eps}_H(k)\leq a_n\ll b_n\leq F^C_H(k,n) $ whp, and
\begin{equation}
F^{2C}_H(n)
\leq 
\frac{1+\eps}{1-\eps}F^C_H(k,n)
\leq
\frac{1+\eps}{(1-\eps)^2}\tilde F_H(n),
\end{equation}
with high probability, which gives the first inequality of (\ref{ineqs:sandwich}).

But since (\ref{ineqs:sandwich}) is valid for any CDF $\tilde G$ with $\tilde G(x) = x+o(x)$ as $x\to 0$, in particular it is valid for $G$, and thus $1-4\eps\leq \tilde F_H(n)/F^{2C}_H(n) \leq 1+\eps$ as well.
It follows that $1-6\eps \leq \tilde F_H(n)/F_H(n)\leq 1+6\eps$ whp. Since $\eps$ was arbitrary, $\tilde F_H(n)/F_H(n)\to 1$ in probability. \end{proof}

\section{Examples of unbalanced cover}
\label{section:unconcentrated-cover}
We'll conclude with examples of cover problems where the upper and lower bounds on $C_H$ don't match, and where $C_H$ is not sharply concentrated. Recall that the lower bound on $C_H$ was of order $n^{1-1/\max(d_H,\Delta)}$, while for factors it was $n^{1-1/d^*}$ (with $d_H, \Delta\leq d^*$).

Why are the lower bounds for factors and covers different? 
If $d_H<d^*$, then $H$ has a denser subgraph $H^*$, and the minimal $H$-cover might have many copies of $H$ overlapping in the same copy of $H^*$.
In an $H$-factor there are at least $\Omega(n)$ vertices lying in some copy of $H^*$ (because $t$ disjoint copies of $H$ contain at least $t$ disjoint copies of $H^*$), while in an $H$-cover only $1$ copy of $H^*$ is guaranteed. 

For the sake of simplicity, let's compare with the threshold for the appearance of an $H$-cover in $G_{n,p}$. The threshold for the existence of a collection of copies of $H^*$ that cover at least $\Omega(n)$ vertices is $p=n^{-\beta}$ with $\beta= 1/d^*=\min_{H'\subseteq H} \frac{v_H-1}{e_H}$. But the threshold for the appearance of at least one copy of $H^*$ is lower, with  $\beta=\min_{H'\subseteq H} \frac{v_H}{e_H}$.

For example, consider $H=K_4 +K_2$ (disjoint union of the complete graph on $4$ vertices and an edge). Here the $1$-density of $H$ is $d_H=1.4$, while the maximum $1$-density of a subgraph is $d^*=2$ (the $K_4$). The maximum $0$-density is $1.5$ (again, the $K_4$). So $\max(d_H,\Delta)=1.5$, and
\cref{lowerbound:cover} gives the lower bound $C_H(n)= \Omp(n^{1/3})$, while \cref{upperbound-Hfactor} gives the upper bound $C_H(n)= \Op(n^{1/2})$.

For this $H$, the lower bound is tight: the cheapest $H$-factor will typically be the cheapest $K_4$
together with the cheapest cover of the remaining $n-4$ vertices by edges.
Define the random variable $Z$ by the lowest weight of a copy of $K_4$ in $K_n$. With a first and second moment method counting the number of $K_4$'s cheaper than $cn^{-2/3}$, one can show that $Z=\Tp(n^{-2/3})$, but $\Pr(Z\leq cn^{-2/3})$ is bounded away from both $0$ and $1$ for any $c$. In other words, $Z$ is not sharply concentrated.

A red-green split argument like in \cref{redgreensplit} with $t=1/2$ leads to a coupling such that $C_H\leq (n-4) Z+2C_{K_2}$ (where $C_H=C_H(n)$ and $C_{K_2}=C_{K_2}(n-4)$), because the smallest number of copies of $H$ that can overlap in the same copy of $K_4$ while also covering all $n$ vertices is $(n-4)/2$. And $C_{K_2}(n-4)=\Op(1)$, which can be seen by either considering a greedy algorithm or using \cite{matching-zeta-2}.
On the other hand, $C_H(n)\geq nZ/6$, because any cover contains at least $n/6$ copies of $H$ that each must contain a $K_4$, and each such copy has weight at least $Z$.
Together this gives us that $ \frac{1}{6}\leq C_H(n)/nZ\leq 1+\op(1)$ whp. Hence $C_H(n)=\Tp(n^{1/3})$, but since $Z$ is not sharply concentrated, neither is $C_H$.

One might guess that this pathological behavior is due to $H$ being disconnected, but it occurs even for some connected graphs. For instance, if $H$ is a $(5,2)$-lollipop graph: a complete graph $K_5$, with a path $P_2$ away from one of the vertices of the clique. There are $7$ vertices and $12$ edges, so $d_H=2$. Since the densest subgraph is the $K_5$,  $\Delta=2$ and $d^*=5/2$. From \cref{mainthm:cover}, the asymptotics of $C_H$ is then between $n^{0.5}$ and $n^{0.6}$.
Here a near-optimal $H$-cover can be found that is a single $K_5$ together with a large collection of paths from this clique.
Consider $G_{n,p}$ with $p=n^{-1/2+\eps}$ for some small $\eps>0$. With high probability, this graph contains a $K_5$ and has diameter $2$. 
Hence $C_H=\Op(np)=\Op(n^{1/2+\eps})$, which is arbitrarily close to the lower bound from \cref{mainthm:cover}.

\begin{small}
\bibliographystyle{abbrv}
\bibliography{FactorBib}

\begin{thebibliography}{10}

\bibitem{matching-zeta-2}
D.~J. Aldous.
\newblock The {$\zeta(2)$} limit in the random assignment problem.
\newblock {\em Random Structures Algorithms}, 18(4):381--418, 2001.

\bibitem{probmeth}
N.~Alon and J.~H. Spencer.
\newblock {\em The probabilistic method}.
\newblock John Wiley \& Sons, 2016.

\bibitem{threshold:connectivity}
P.~Erd\H{o}s and A.~R\'{e}nyi.
\newblock On random graphs. {I}.
\newblock {\em Publ. Math. Debrecen}, 6:290--297, 1959.

\bibitem{threshold:matching}
P.~Erd\H{o}s and A.~R\'{e}nyi.
\newblock On the existence of a factor of degree one of a connected random
  graph.
\newblock {\em Acta Math. Acad. Sci. Hungar.}, 17:359--368, 1966.

\bibitem{tree-zeta-3}
A.~M. Frieze.
\newblock On the value of a random minimum spanning tree problem.
\newblock {\em Discrete Appl. Math.}, 10(1):47--56, 1985.

\bibitem{jsl-randomgraphs}
S.~Janson, T.~{\L}uczak, and A.~Ruci\'{n}ski.
\newblock {\em Random graphs}.
\newblock Wiley-Interscience Series in Discrete Mathematics and Optimization.
  Wiley-Interscience, New York, 2000.

\bibitem{JohKahVu08}
A.~Johansson, J.~Kahn, and V.~Vu.
\newblock Factors in random graphs.
\newblock {\em Random Structures Algorithms}, 33(1):1--28, 2008.

\bibitem{KOMLOS198355}
J.~Komlós and E.~Szemerédi.
\newblock Limit distribution for the existence of hamiltonian cycles in a
  random graph.
\newblock {\em Discrete Mathematics}, 43(1):55--63, 1983.

\bibitem{threshold:hamiltoncycle1}
A.~Kor{\v{s}}unov.
\newblock Solution of a problem of p. erdos and a. r{\'e}nyi on hamiltonian
  cycles in undirected graphs.
\newblock In {\em Dokl. Akad. Nauk SSSR}, volume 228, pages 529--532, 1976.

\bibitem{threshold:hamiltoncycle2}
L.~P\'{o}sa.
\newblock Hamiltonian circuits in random graphs.
\newblock {\em Discrete Math.}, 14(4):359--364, 1976.

\bibitem{rucinski-partialfactor}
A.~Ruci\'{n}ski.
\newblock Matching and covering the vertices of a random graph by copies of a
  given graph.
\newblock {\em Discrete Math.}, 105(1-3):185--197, 1992.

\bibitem{talagrand1995concentration}
M.~Talagrand.
\newblock Concentration of measure and isoperimetric inequalities in product
  spaces.
\newblock {\em Publications Math{\'e}matiques de l'Institut des Hautes Etudes
  Scientifiques}, 81:73--205, 1995.

\bibitem{walkup1979expected}
D.~W. Walkup.
\newblock On the expected value of a random assignment problem.
\newblock {\em SIAM Journal on Computing}, 8(3):440--442, 1979.

\bibitem{tsp-204}
J.~W\"{a}stlund.
\newblock The mean field traveling salesman and related problems.
\newblock {\em Acta Math.}, 204(1):91--150, 2010.

\end{thebibliography}

\end{small}

\end{document}